\documentclass[12pt]{amsart}
\usepackage{amsmath}
\usepackage{amsfonts}
\usepackage{amsthm}
\usepackage{amssymb}
\usepackage[english]{babel}
\usepackage[ansinew]{inputenc}
\usepackage[all]{xy}
\usepackage{color}

\theoremstyle{definition}     

\newtheorem{defi}{Definition}[section]
\newtheorem{remark}[defi]{Remark}

\newtheorem{example}[defi]{Example}
\newtheorem{definition}[defi]{Definition}

\theoremstyle{plain}

\newtheorem{theorem}[defi]{Theorem}
\newtheorem{corollary}[defi]{Corollary}
\newtheorem{lemma}[defi]{Lemma}
\newtheorem{proposition}[defi]{Proposition}

\input xy
\xyoption{all}

\title{The pseudo-fundamental group scheme}

\author{Marco Antei}\thanks{Marco Antei thanks the project TOFIGROU (ANR-13-PDOC-0015-01)}

\address{Universidad de Costa Rica, Ciudad universitaria Rodrigo Facio Brenes, Costa Rica.}

\email{marco.antei@ucr.ac.cr}

\author{Arijit Dey}

\address{Indian Institute of Technology Madras, Sardar Patel Road, Chennai- 600036, India.}

\email{arijit@iitm.ac.in}

\date{}
\begin{document}

\begin{abstract}
Let $X$ be any scheme defined over a Dedekind scheme $S$ with a given section $x\in X(S)$. We prove the existence of a pro-finite $S$-group scheme $\aleph(X,x)$ and a universal $\aleph(X,x)$-torsor dominating all the pro-finite pointed torsors over $X$. Though $\aleph(X,x)$ may not be unique in general it still can provide useful information in order to better understand $X$. In a similar way we prove the existence of a pro-algebraic $S$-group scheme $\aleph^{\rm alg}(X,x)$ and a $\aleph^{\rm alg}(X,x)$-torsor dominating all the pro-algebraic and affine pointed torsors over $X$. The case where $X\to S$ has no sections is also considered. 
\end{abstract}
\maketitle


\tableofcontents
\bigskip

\section{Introduction}
\label{sez:Intro}
The existence of a group scheme classifying all finite torsors over a given scheme $X$ was first  conjectured by Grothendieck in \cite[Chapitre X]{SGA1}. It was first Nori who proved it (cf.  \cite{Nor} and \cite{Nor2}) and called it the \emph{fundamental group scheme} (FGS). Infact in his thesis \cite{Nor2}, Nori gave two possible ways to construct the FGS. In the first method, he constructed it for a connected, proper and reduced scheme $X$ defined over a perfect\footnote{In \cite{AnEm} it has been pointed out that the perfectness assumption for the field $k$ was only needed to ensure that $H^0(X,\mathcal{O}_X)=k$, so instead of considering only perfect fields one can take \emph{any} field $k$ with the additional assumption on the scheme $X$ that $H^0(X,\mathcal{O}_X)=k$} field $k$ equipped with a section $x\in X(k)$,  as a $k$-group scheme $\pi(X,x)$ naturally associated to the neutral tannakian category of essentially finite vector bundles over $X$. In the second method he assumes $X$ to be reduced and connected but not necessarily proper, with arbitrary underlying field $k$, such that $X$ has at least one $k$-rational point $x$. With these assumptions he proves the existence of  $\pi(X,x)$ by showing that the category of all finite pointed torsors over $X$ is cofiltered which is a necessary and sufficient condition for existence of FGS \cite[Proposition 2]{Nor2}. Clearly his second method works in much more general set-up than the first one. The second method has also been generalized for schemes defined over Dedekind schemes (see \cite{GAS}, \cite{AEG}). More precisely the main result in \cite{AEG} is the existence of $\pi(X,x)$ when $X\to S$ is separated, faithfully flat, of finite type and either for all $s\in S, X_s$ is reduced, or for all $x\in X\backslash X_{\eta}$, $\mathcal{O}_x$ is integrally closed (here $\eta$ denotes the generic point). Under extra assumptions a quasi-finite version $\pi^{\text{qf}}(X,x)$ has also been studied, i.e. a fundamental group scheme classifying all the quasi-finite pointed torsors over $X$.    

After Nori's work the next step was to carry forward FGS construction for non-reduced pointed schemes defined over a field or more generally over a Dedekind scheme $S$ with a section $x \in X(S)$. It is known that the category of all finite torsors over a non-reduced scheme $X$ defined over $S$ is not cofiltered in general (\cite[Reamrk 1.3. (iii)]{LZ}). In \cite{Antei2}, the first author made an attempt to construct FGS in this setting by showing that the category of finite pointed \emph{Galois} torsors (see Definition \ref{defGal}) is cofiltered. Unfortunately proof in \cite{Antei2} contains a mistake, in fact here we give an actual counterexample to his claim (cf. Example \ref{exPoint}). In this paper we keep the same idea of \cite{Antei2} i.e. considering Galois torsors instead of taking all of them but in a larger category of pro-finite torsors.  
Only in this new environment we are able to show that there exists a Galois torsor (which will be called \emph{universal}) dominating all the finite Galois torsors (cf. Theorem \ref{teoPseudo1}). The structural group scheme of this universal torsor will be denoted by $\aleph(X,x)$ and called the \emph{pseudo-fundamental group scheme} (PFGS) of $X$ at the point $x$. In general this new object PFGS need not be unique though any two such PFGS are dominated by a third one.  
However, in \cite{BV} Borne and Vistoli generalized the notion of Nori's FGS to fundamental gerbe, which applies to schemes, algebraic stacks, and more generally to a fibered category  even in absense of rational points.  Using Nori's approach they proved that a fibred category $\mathcal X$ has a fundamental gerbe if and only if it is inflexible \cite[Theorem 7. and 7.13]{BV}.

It is clear that whenever $X$ admits a fundamental group scheme $\pi(X,x)$ then it coincides with $\aleph(X,x)$, up to a unique isomorphism. 
It is natural to wonder how the natural morphism $\pi(X_{\rm red},x)\to \aleph(X,x)$ behaves, provided $\pi(X_{\rm red},x)$ exists. In the \'etale case this morphism is known to be an isomorphism (cf. \cite[I, Th\'eor\`eme 8.3]{SGA1}).  We will discuss this at the end of \S \ref{sez:withPoint}. 

In \cite{BV2}, Borne and Vistoli proved that for a  fibred category $\mathcal X$, under some mild assumptions, has a virtually abelian and a virtually unipotent fundamental group gerbe. Those are pro-algebraic (not necessarily pro-finite) group gerbes. 
The techniques used in this paper will also work in this general setting. In particular this allows us to construct a universal torsor dominating \emph{all} the Galois objects in the category of pointed algebraic torsors (i.e. those torsors whose structural group scheme is affine and of finite type), implying the existence of a \emph{algebraic pseudo-fundamental group scheme} $\aleph(X,x)^{\rm alg}$ satisfying similar properties. This is the biggest possible affine torsor over $X$, thus dominating all the others constructed so far. The importance of introducing $\aleph(X,x)^{\rm alg}$ is also the possibility to have a finer invariant than $\pi(X,x)$; indeed in Lemma \ref{lemFiner} we show that if $X$ is a smooth and connected projective scheme over a field $k$, then $\aleph(X,x)^{\rm alg}$ is trivial if and only if $X$ is a point. It may thus be useful to study this object from an anabelian point of view. Moreover, as pointed out in \cite{BV2}, an algebraic (not pseudo) fundamental group scheme does not exist in general: this means that there is no (unique) universal torsor dominating all the algebraic affine torsors, but if we are willing to give up on the unicity of our (uni)versal object, then we can find a \emph{biggest} torsor dominating all the others.  Furthermore, again, the same techniques are used in \S \ref{sez:noPoint} to define a \emph{non pointed} version of pseudo-fundamental group schemes in the pro-finite and the pro-algebraic environments. 
In the classical case Nori used two properties very crucially to prove existence and uniqueness of the fundamental group scheme $\pi(X,x)$. One is underlying scheme $X$ is reduced and other is torsors are finite and pointed. This helped him to prove that the category of finite torsors is cofiltered, hence the existence of the $\pi(X,x)$-universal torsor \cite [Chapter II, Proposition 2]{Nor}. Whereas our approach allows us to drop these two strong conditions and we show the existence of  $\aleph(X)$ which classifies all pro-finite torsors having the desired universal property. The generality of our method also allows us to construct such object 
$\aleph(X)^{alg}$ in pro-algebraic setting. This last construction also represents an alternative to the \emph{fundamental groupoid schemes} and \emph{fundamental gerbes} already considered in \cite{EH} and \cite{BV} respectively when we want to bypass the existence of a rational point.

\section{Preliminaries}
\label{sez:Revenge}

Let $S$ be any Dedekind scheme (e.g. the spectrum of a field or a discrete valuation ring) and $\eta=Spec(K)$ be its generic point. Let $X$ be a scheme over $S$ endowed with a $S$-valued point $x:Spec(S)\to X$. 
A triple $(Y,G,y)$ over $X$ is a \textit{fpqc}-torsor $Y\to X$, under the (right) action of a flat affine $S$-group scheme $G$ together 
with a $S$-valued point $y\in Y_x(S)$. The morphism between two triples $(Y,G,y)\to(Y',G',y')$ are morphisms of $S$-schemes 
$\alpha:G\to G'$ and a $G$-equivariant morphism $\beta:Y\to Y'$ such that $\beta(y)=y'$.
The category whose objects are triples $(Y,G,y)$ with the additional assumption that $G$ is finite and flat is denoted by $\mathcal{P}(X)$. 
We denote by $Pro-\mathcal{P}(X)$ the pro-category of $\mathcal{P}(X)$ whose objects are projective limits of objects in $\mathcal{P}(X)$ and as usual for any two objects 
$(Y,G,y)=\varprojlim_i (Y_i,G_i,y_i)$, $(T,M,t)=\varprojlim_j (T_j,M_j,t_j)$ in $Pro-\mathcal{P}(X))$, morphisms between them is given by 
$$Hom\left((Y,G,y), (T,M,t)\right)\,=\,\varprojlim_j\varinjlim_i Hom_{\mathcal{P}(X)}\left((Y_i,G_i,y_i),(T_j,M_j,t_j)\right).$$

$Pro-\mathcal{P}(X)$ is a full subcategory of the category $\overline{\mathcal{P}(X)}$ whose objects are triples of torsors under the action of affine and flat group schemes. The same is true for $Pro-\mathcal{P}^{\text{alg}}(X)$, where $\mathcal{P}^{\text{alg}}(X)$ is the category of triples $(Y,G,y)$, as before, with the only difference that $G$ is now flat and affine. 

\begin{definition}\label{defGal} We say that an object $(Y,G,y)$ of $\mathcal{P}(X)$ (resp. of $Pro-\mathcal{P}(X)$)  over $X$ is Galois relatively to $\mathcal{P}(X)$ (resp. to $Pro-\mathcal{P}(X)$) if for every triple $(Y',G',y')$ of $\mathcal{P}(X)$ (resp. of $Pro-\mathcal{P}(X)$) and every morphism $(Y',G',y')\to (Y,G,y)$ the group scheme morphism $G'\to G$ is faithfully flat (or, equivalently the morphism $Y'\to Y$ is faithfully flat). The full subcategory of $\mathcal{P}(X)$ (resp. of $Pro-\mathcal{P}(X)$) whose objects are Galois triples is denoted by $\mathcal{G}(X)$ (resp. $\mathcal{G}(X)'$). In a similar way we define $\mathcal G^{\text{alg}}(X)$ and $\mathcal G^{\text{alg}}(X)'$ as full subcategories of $\mathcal{P}(X)$ and $Pro-\mathcal{P}(X)$ respectively. 
\end{definition}

Notice that a Galois triple is \emph{minimal} in the sense that any morphism $(Z,H,z)\to (Y,G,y)$ where $(Y,G,y)$ is Galois and $H\to G$ is generically a closed immersion is necessarily an isomorphism. Full details of its construction will be given in Proposition \ref{propOldRed2}.

\begin{remark}\label{remGalProj} A projective limit of objects of $\mathcal{G}(X)$ denoted as $Pro-\mathcal{G}(X)$ 
is an object of $\mathcal{G}(X)'$ [\cite{EGAIV-3}, \text{section} 8.3.8]. We can restate this saying that $Pro-\mathcal{G}(X)$ is a full 
subcategory of $\mathcal{G}(X)'$. It is not clear to us if the inclusion functor is essentially surjective. 
\end{remark}

If the category $\mathcal{G}(X)$ was cofiltered we could easily deduce the existence of a universal torsor projective limit of all the objects in $\mathcal{G}(X)$ (unique up to a unique isomorphism). Unfortunately it is not true in general when $X$ is not reduced. Indeed we provide an example where an object of $\mathcal{G}(X)$ where $X=Spec(k[x]/x^2)$ has a non trivial automorphism, which implies that $\mathcal{G}(X)$ is not cofiltered:  

\begin{example}
\label{exPoint}

Here we show that if $X=Spec(k[x]/x^2)$, where $k$ is a field of characteristic $2$, the category $\mathcal{G}(X)$ is not cofiltered. It is sufficient to find a $k$-group scheme $G$ and a pointed $G$-torsor $Y$ in $\mathcal{G}(X)$ and an automorphism (in $\mathcal{G}(X)$) different from the identity. We choose $G:=\alpha_2=Spec(k[x]/x^2)$ and $Y:=Spec(k[x,y]/(x^2,y^2+x))$ is a $G$-torsor pointed in the origin. It is not trivial as for all $a\in k[x]/x^2$, $x\neq a^2$. A non-trivial pointed $\alpha_2$ torsor is therefore necessarily Galois. The right action of $G$ on $Y$ can be described as a coaction as follows:

$$\begin{array}{ccc}
\rho: k[x,y]/(x^2,y^2+x) & \to & k[x]/x^2 \otimes_k k[x,y]/(x^2,y^2+x)\\
x & \mapsto & 1\otimes x  \\
y & \mapsto & x\otimes 1 + 1 \otimes y   
\end{array}$$

and this action is giving the following isomorphism making $Y$ a $G$-torsor:

$$\begin{array}{ccc}
\frac{k[x,y]}{x^2,y^2+x}\otimes_{k[x]/x^2} \frac{k[x,y]}{x^2,y^2+x} & \to & \frac{k[x]}{x^2} \otimes_k \frac{k[x,y]}{x^2,y^2+x}\\
1\otimes x & \mapsto & 1\otimes x  \\
1\otimes y & \mapsto & 1 \otimes y   \\
y \otimes 1 & \mapsto & x\otimes 1 + 1 \otimes y.   
\end{array}$$

Now we consider the morphism of $k[x]/x^2$-algebras
$$\varphi^{\#}: k[x,y]/(x^2,y^2+x)\to k[x,y]/(x^2,y^2+x), \qquad  y\mapsto x+y$$
and we observe that it commutes with the coaction as $(id\otimes \varphi^{\#}) \rho = \rho \varphi^{\#} $. The induced morphism of $X$-schemes $\varphi:Y\to Y$ and the identity morphism on $G$ give a morphism in  $\mathcal{G}(X)$, different from the identity. Hence $\mathcal{G}(X)$ is not cofiltered. 
\end{example}
\begin{remark}
The existence of the above counter-example can be also seen as a consequence of a more geenral fact: If $Y\rightarrow X$ is a $G$-torsor with $G$ being commutative and connected then any $g\in G(X)\backslash 1_{G(X)}$ is a $G$-equivariant automorphism of $Y$ over $X$ which sends the (unique) $k$-point maps to itself. We choose to give the above constructional proof because of it's simplicity.
\end{remark} 


We recall the following well known definition from category theory which will be used crucially in this paper. 
\begin{definition}
A skeleton of a category $\mathcal {C}$ is a full subcategory $\text{Sk}(\mathcal {C})$ in which every object in $\mathcal {C}$ is isomorphic to an 
object in $\text{Sk}(\mathcal {C})$ and no distinct objects in $\text{Sk}(\mathcal {C})$ are isomorphic in $\mathcal {C}$. 
\end{definition}

\begin{theorem}
A skeleton of $\mathcal {C}$ always exists. Every skeleton of $\mathcal {C}$ is equivalent to $\mathcal {C}$.
\end{theorem}
\proof
Cf. \cite{AHS}, Proposition 4.14.
\endproof

\begin{definition}
A $S$-morphism $i:Y\to Z$ is said to be a generically closed immersion if its restriction $i_{\eta}:Y_{\eta}\to Z_{\eta}$ to $\eta$ is a closed immersion. When $S$ is the spectrum of a field that simply means that $i$ is a closed immersion.
\end{definition}

\begin{remark}\label{remHai}
A $S$-morphism of group schemes $G\to G'$ can be factored into a faithfully flat morphism $G\to Q$, a \emph{model map}\footnote{Notice that a model map is not in general a monomorphism: monomorphisms are stable after base change, and special fibers of model maps are often not monomorphisms.} (i.e. generically an isomorphism, as defined for instance in \cite{BLR}) $Q\to M$ and a closed immersion $M\to G'$. When $S$ has dimension 0 then $Q\to M$ is an isomorphism.
\end{remark}

\begin{proposition}\label{propOldRed2} Given an object $(Y,G,y)$ of $Pro-\mathcal{P}(X)$, there exists an object $(T,H,t)$ of $\mathcal{G}(X)'$ and a
morphism $(T,H,t)\to (Y,G,y)$ where $T\to Y$ (or equivalently $H\to G$) is a generically closed immersion. We say in this case that $(T,H,t)$ is contained in $(Y,G,y)$.
\end{proposition}
\proof

Let $C_Y$ denote the category whose objects are $(A,H,a,f)$ where $(A,H,a)$ is an object in ${Pro}-\mathcal P(X)$
and $f:\, A \longrightarrow Y$ is a generically closed immersion in ${Pro}-\mathcal P(X)$ which takes the point $a$ to $y$. By abuse of notation 
we denote an object in $C_Y$ by $(A,f)$. We first need to prove that isomorphism classes of objects in $C_Y$ forms a set. First we claim that the possible such $H$'s forms a set, indeed the coordinate ring of $H$ are sub-algebras of quotients
of the coordinate ring of generic fiber of $G$. Given such a $H\to G$ we have to deal with the possible $H$-torsors $Q$ with a map $Q\to  Y$. Choose a presentation $H=\varprojlim_q
H_q$. From the map $H\to G=\varprojlim_i G_i$ we find indexes $q_i$ and maps $H_{q_i}\to G_i$. The $H_{q_i}$-torsor induced by the $Q_i$'s of finite type over $X$ and as a map of the $G_i$-torsor induced by $Y$.
Since the set of $q$'s is given we have a set of possible $Q$. Since sections $Q(S)$ are a
set, we are done.

Now we need to prove that between two objects $(A,f)$ and 
$(B,g)$ there is at most one morphism: indeed if such a morphism exists then it turns out to be generically closed immersions $h:A \longrightarrow B$ such that following diagram commutes:
 $$
 \xymatrix{
 A  \ar[d]_{h} \ar[r]^f &Y \\ 
 B  \ar[ru]_g & 
 }
 $$
 Note that the only endomorphism of any object in $C_Y$ is the identity: it is easy to verify over the generic point $\eta$ of $S$ then we observe that since $A$ is flat, then the same holds globally. Then the isomorphism classes of objects of $C_Y$ form a partially ordered set. It is now an easy application of Zorn's Lemma the existence of a minimal element $(Y_{\text{min}}, G_{\text{min}}, y_{\text{min}})$ in $Pro-\mathcal{P}(X)$. To show that it is Galois, assume it is not, then there exists (as recalled in Remark \ref{remHai}) a triple $(U,M,u)$ in $Pro-\mathcal{P}(X)$ and a morphism $(U,M,u)\to (Y_{min},G_{min},y_{min})$ which is not faithfully flat. 
Hence it will factor through a faithfully flat morphism $(U,M,u)\to (U',M',u')$ and a generically closed immersion $(U',M',u')\to (Y_{min},G_{min},y_{min})$ 
contradicting the minimality, so we can set $(T,H,t)\,:=\,(Y_{min},G_{min},y_{min})$).
\endproof

\begin{corollary}\label{corOldRed2} Given two objects $(Y_i,G_i,y_i)$ in $\mathcal{G}(X)'$, $i=1, 2$, there exists an object 
$(Y_3,G_3,y_3)\,\in\,\mathcal{G}(X)'$ 
with morphisms $(Y_3,G_3,y_3)\to (Y_i,G_i,y_i)$, for $i=1, 2$.
\end{corollary}
\proof
It is sufficient to consider $(Y_1\times_X Y_2,G_1\times_S G_2,y_1\times_k y_2)$ and then to take a Galois triple contained in it following 
Proposition \ref{propOldRed2}.
\endproof

\begin{remark}
Note that in the proof of Proposition \ref{propOldRed2} we have not used any property of finite group schemes. Therefore instead of working with $\mathcal{P}(X)$ we can work with $\mathcal{P}^{\text{alg}}(X)$ and the proof of Proposition \ref{propOldRed2} goes through exactly in similar fashion, hence the 
Corollary \ref{corOldRed2}. 
\end{remark}

\section{The pseudo-fundamental group scheme }
\label{sez:new}
We first study the problem of the existence of a group scheme classifying all finite (resp. affine and of finite type)  torsors in the ``classical'' case of pointed schemes, considering maps between torsors sending the marked point of the source torsor to the marked point of the target. For this reason at the end of \S \ref{sez:withPoint} we will be able to compare $\aleph(X,x)$ to $\pi(X_{\rm red},x)$, the latter being Nori's fundamental group scheme. In \S \ref{sez:noPoint} we will provide a short overview on the category of non pointed torsors.

\subsection{The case of pointed torsors}
\label{sez:withPoint}

\begin{definition}\label{defPseudoFirst}
The $S$ scheme $X$ has a \emph{pseudo-fundamental group scheme} (PFGS) $\aleph(X,x)$ if there is a triple $(\widehat{X},\aleph(X,x),\widehat{x})$ in the category $\mathcal{G}(X)'$ 
such that for each object $(Y,G,y)$ in $\mathcal{G}(X)$ there is a morphism
$(\widehat{X},\aleph(X,x),\widehat{x})\to (Y,G,y)$. In this case $\widehat{X}$ is called 
the universal $\aleph(X,x)$-torsor over $X$ pointed in $\widehat{x}$.
\end{definition}


\begin{remark}\label{propPseudoFirst}
Though the PSGS may not be unique whenever $(\widehat{X},\aleph(X,x),\widehat{x})$ and 
$(\widehat{X}',\aleph(X,x)',\widehat{x}')$ are two PFGS triples for $X$ then there exists a third one $(\widehat{X}'',\aleph(X,x)'',\widehat{x}'')$ dominating both. This is an easy consequence of Corollary \ref{corOldRed2}. However this does not imply the existense of a (even) bigger one dominating all of them. Moreover if the PFGS of $X$ is known to be finite then all the universal triples are of course all isomorphic (but the isomorphism may not be unique unlike in Nori's case.)
For finite type torsors one can similarly define the algebraic pseudo-fundamental group scheme (APFGS) $\aleph^{\text{alg}}(X,x)$ as an object 
$(\widehat{X}^{\text{alg}},\aleph^{\text{alg}}(X,x),\widehat{x})$ in ${\mathcal{G}}^{\text{alg}}(X)'$
\end{remark}

\begin{definition} For two triples $(Y_1,G_1,y_1)$ and $(Y_2,G_2,y_2)$ in $\mathcal{G}(X)'$ (resp. $\mathcal{G}^{\text{alg}}(X)'$) we say that $(Y_1,G_1,y_1)$ dominates $(Y_2,G_2,y_2)$ if there 
exists a (maybe not unique) morphism $(Y_1,G_1,y_1)\to (Y_2,G_2,y_2)$.
\end{definition}

\begin{theorem}\label{teoPseudo1}
Let $X$ be a scheme over $S$ with a $S$-valued point $x$. Then $X$ has a PFGS (resp. APFGS)  $\aleph(X,x)$ (resp. $\aleph^{\text{alg}}
(X,x))$).
\end{theorem}
\proof

The proof for the existence of PFGS and APFGS are exactly similar. Here we give a proof for the existence of PFGS. In order to simplify the discussion we only consider finite and pro-finite torsors, discussion on (pro-)algebraic torsors will be identical (mutatis mutandis). Moreover a triple $(Z,G,y)$ (a $G$-torsor $Z$ pointed in $y$ over $x$) will simply be denoted by $Z$. Let us consider $Sk(q\mathcal{G}(X))$ the skeleton of the quotient category (with one single morphism) of $\mathcal{G}(X)$, the latter being the category of Galois finite torsors. We put on $I:=Ob(Sk(q\mathcal{G}(X)))$ (which is a set) a well order $(Ob(Sk(q\mathcal{G}(X))),\prec)$: 

$$Z_1, Z_2, ... , Z_{\alpha}, ...$$ 

then we argue like this: we call $Z'_2$ a torsor in $\mathcal{G}'(X)$ dominating both $Z_1$ and $Z_2$ (it always exists by Corollary \ref{corOldRed2}). Then we call $Z'_3$ a torsor in $\mathcal{G}'(X)$ dominating both $Z'_2$ and $Z_3$ and we go on like this thus obtaining a chain in $\mathcal{G}'(X)$. Applying the functor $\text{Sk} \circ q$ we obtain a chain of morphisms 

$$... \to  Z'_{\beta}\to ... \to Z'_{\alpha}\to ... \to  Z'_3\to Z'_2\to Z_1 (\dagger)$$
in $Sk(q(\mathcal{G}'(X)))$. Without loss of generality we can assume that $Z_i'$'s are distinct. This because in $Sk(q(\mathcal{G}'(X))$, for any two object $A$ and $B$, $A = B$ if and only if $Hom(A,B)$ and $Hom(B,A)$ are both nonempty. Lifting the chain $\dagger$ to any chain in $\mathcal{G}'(X)$ we compute the projective limit $\hat Z$. Since  $Pro(Pro(P(X)))$ is equivalent to $Pro(P(X))$ (\cite{DDH}, Theorem 2.17) 
$\hat Z$ is an element 
of $Pro(P(X))$ and by applying Proposition \ref{propOldRed2} we can assume it to be an element of $\mathcal{G}(X)'$. 
\endproof






Let now $X$ be a connected $S$-scheme of finite type with a given section $x\in X(S)$. Let $X_{\rm red}$ be its reduced part. As precised in \S \ref{sez:Intro} we assume that for $X_{\rm red}$ we are able to build the fundamental group scheme $\pi(X_{\rm red},x)$; this is always possible when $dim(S)=0$. We choose a PFGS $\aleph(X,x)$ and a universal $\aleph(X,x)$-torsor $\widehat{X}\to X$, pointed in $\widehat{x}\in \widehat{X}_x(k)$. We consider its pullback $\widehat{X}_X$ over $X$ and the unique morphism of torsors $$\varphi_{\rm red}:X^N\to \widehat{X}$$ where $X^N\to X$ is the (``N'' stands for Nori) universal $\pi(X_{\rm red},x)$-fundamental group scheme. Though in characteristic $0$ this morphism is known to be an isomorphism, in positive characteristic this is no longer true: for instance when $S=Spec(k)$, $k$ being a field with $char(k)=2$, $X=Spec(k[x]/x^2)$, in Example \ref{exPoint} we recalled that over $X=Spec(k[x]/x^2)$ there are non trivial Galois pointed torsors while over $X_{\rm red}=Spec(k)$ there are only trivial pointed torsors. So in this case $\varphi_{\rm red}$ is trivially a closed immersion. In general we have the following:
\begin{proposition}Let $X$ be any affine scheme of finite type over a Dedekind scheme $S$, endowed with a section $x\in X(S)$.  and let $X_{\rm red}$ its reduced part, for which we assume it admits a fundamental group scheme $\pi(X_{\rm red},x)$. Let moreover $\aleph(X,x)$ be a PFGS of $X$ then the morphism $i:X_{\rm red}\to X$ induces a closed immersion $$\varphi_{\rm red}:\pi(X_{\rm red},x)\to \aleph(X,x).$$
\end{proposition}
\proof
This follows essentially from \cite[Section 3.2]{Antei2}, indeed we know that every finite torsor over  $X_{\rm red}$ can be extended over $X$, then in particular we get a morphism $q:\aleph(X,x)\to \pi(X_{\rm red},x)$ such that $q\circ \varphi_{\rm red}=id_{\pi(X_{\rm red},x)}$. This implies that $\varphi_{\rm red}$ is a closed immersion.
\endproof
In a similar way one can study the morphisms between $\aleph^{\rm alg}(X_{\rm red},x)$ and $\aleph^{\rm alg}(X,x)$.

\begin{lemma}
\label{lemFiner}
If $X$ is a smooth and connected projective scheme over a field $k$ then $\aleph^{\rm alg}(X,x)$ is trivial if and only if $X=Spec(k)$. 
\end{lemma}
\proof If $X$ is not a point it is sufficient to consider the closed immersion $j:X\to \mathbb{P}^n_k$, some $n\in \mathbb{N}$, and to observe that $j*(\mathcal{O}_{\mathbb{P}^n_k}(1))$ is a non trivial line bundle. This gives rise to a non trivial $\mathbb{G}_m$-torsor $Y$ over $X$. 
\endproof

This is false and well known for $\pi(X,x)$: for instance $\pi(\mathbb{P}^n_k,x)$ is known to be trivial.

\subsection{The case of non pointed torsors}
\label{sez:noPoint}
In \S \ref{sez:Intro} we made clear that we first defined the pseudo fundamental group scheme giving a $S$-valued point $x$ on $X$ in order to compare it to Nori's fundamental group scheme whose constructions (both the tannakian and the pro-finite) always need a given point. However when we work over non algebraically closed fields or Dedekind schemes it can be useful to have a similar object even when such a point does not exist. The reader certainly observed that the proofs of \S \ref{sez:withPoint} still holds if the base scheme $X$ and torsors are not pointed. Without repeating the proofs we only introduce new definitions and recall the main properties following same arguments of \S \ref{sez:withPoint}. Here $\mathcal{T}(X)$ and $Pro-\mathcal{T}(X)$ will denote respectively the category of finite torsors over $X$ and that of pro-finite torsors over $X$.

\begin{definition}\label{defGalglob} We say that an object $(Y,G)$ of  $Pro-\mathcal{T}(X)$ is Galois if for every object $(Y',G')$ of $Pro-\mathcal{T}(X)$ and every morphism $(Y',G')\to (Y,G)$ the group scheme morphism $G'\to G$ is faithfully flat (or, equivalently the morphism $Y'\to Y$ is faithfully flat). The full subcategory of $Pro-\mathcal{T}(X)$ whose objects are Galois is denoted by $\mathcal{F}(X)$.
\end{definition}

\begin{definition}\label{defPseudoSecond}
$X$ has a \emph{global pseudo-fundamental group scheme} $\aleph(X)$ if there is a pair $(\widehat{X},\aleph(X))$ in the category $\mathcal{F}(X)$ 
such that for each object $(Y,G)$ of $\mathcal{F}(X)$ there is a morphism $(\widehat{X},\aleph(X))\to (Y,G)$. In this case $\widehat{X}$ is called 
the universal $\aleph(X)$-torsor over $X$.
\end{definition}
Again it is clear by this definition that whenever $\aleph(X)$ and $\aleph(X)'$ are two distinct global pseudo-fundamental group schemes then we have two (maybe not unique) faithfully flat morphisms $\aleph(X)'\to \aleph(X)$ and $\aleph(X)\to \aleph(X)'$ whose compositions are not necessarily automorphisms.

\begin{theorem}\label{teoPseudoGlob}
Let $X$ be a scheme over a Dedekind scheme $S$. Then $X$ has a global pseudo-fundamental group scheme $\aleph(X)$.
\end{theorem}

In a similar way one can define the \emph{global algebraic pseudo-fundamental group scheme} $\aleph^{\rm alg}(X)$ (with a obvious meaning) of a scheme $X$ without specifying the existence of a section $x\in X(S)$ and verify that the statements just recalled still hold. The following remark will conclude the paper:

\begin{remark}When $S$ is the spectrum of an algebraically closed field, $x\in X$ any point and we assume that $X$ has a fundamental group scheme $\pi(X,x)$ then it is not difficult to prove that $\pi(X,x)$ and $\aleph(X)$ are isomorphic.
\end{remark}

\section*{Acknowledgement} We are grateful to Fabio Tonini for his useful comments, suggestion and a very careful reading which has substantially improved the paper. We would also like to thank Dajano Tossici for many useful discussions. The second-named author wishes to thank Universit\'e de Nice Sophia-Antipolis for providing local hospitality and support. He also wishes to thank Department of Science and Technology, India for providing him the research grant MATRIX.

%

\end{document}